\documentclass[12pt]{article}

\usepackage{amsmath,amsxtra,amssymb,latexsym, amscd}
\usepackage[mathscr]{eucal}

\begin{document}                        
                        \def\DATE{9/13/'03}
\title{ {\huge\bf Spectral measure of Laplacean operators in Paley-Wiener space} }
\def\1{\rule{0cm}{0cm}} \def\qd{\rule{3mm}{3mm}} \def\BB{$\bullet$}
\renewcommand{\arraystretch}{1.25}
\renewcommand{\theequation}{\thesectn.\arabic{equation}}
\def\sce{\setcounter{equation}{0}}  \newcounter{sectn} \newcounter{sbsect}
\def\sect#1{\addtocounter{section}{1}\sce\setcounter{sbsect}{0}%
        \renewcommand{\thesectn}{\thesection}\1\smallskip\\
        {\1\hspace{-2em}\large\bf\thesectn.\qquad #1\smallskip\par}}
\def\subsect#1{\addtocounter{sbsect}{1}\sce%
        \renewcommand{\thesectn}{\thesection:\Alph{sbsect}}\1\smallskip\\
        {\bf\1\hspace{-1.5em}\thesectn.\qquad #1\smallskip\par}}
\newtheorem{Theorem}{THEOREM} \newtheorem{Lemma}[Theorem]{LEMMA}
\newtheorem{Corollary}[Theorem]{COROLLARY}
\def\thm#1#2{\be{Theorem}{\lb{#1} #2}} \def\LEM#1#2{\BE{Lemma}{\LB{#1} #2}}
\def\COR#1#2{\BE{Corollary}{\LB{#1} #2}}
\def\proof{\bigskip\noindent {\sc Proof:}\qquad}
\def\REM{\1\smallskip\par\noindent{\bf REMARK:}\qquad }
\def\qed{\hfill$\quad$\qd\medskip\\} \def\ds{\displaystyle}
\def\LB#1{\label{#1}} \def\BE#1#2{\begin{#1} #2 \end{#1}}
\def\EQ#1#2{\BE{equation}{\LB{#1} #2}} \def\ARR#1#2{\BE{array}{{#1} #2}}
\def\DES#1{\BE{description}{#1}} \def\QT#1{\BE{quote}{#1}}
\def\ENUM#1{\BE{enumerate}{#1}} \def\ITM#1{\BE{itemize}{#1}}
\def\COM#1{\par\noindent{\bf COMMENT:\quad\sl #1}\par\noindent}
\def\mapsfrom{\hbox{$\;{\leftarrow}\kern-.15em{\mapstochar}\:\:$}}
\def\vv{\kern.344em{\rule[.18ex]{.075em}{1.32ex}}\kern-.344em}
\def\RE{\mbox{\rm I\kern-.21em R}} \def\CX{\mbox{\rm \vv C}}
\def\imp{\Rightarrow} \def\emb{\hookrightarrow} \def\wk{\rightharpoonup}
\def\rd{\dot{\1}} \def\d{\cdot} \def\+{\oplus} \def\x{\times}
\def\<{\langle} \def\>{\rangle} \def\o{\circ} \def\at#1{\Bigr|_{#1}}
\def\cd{\partial} \def\grad{\nabla} \def\L{\left} \def\R{\right}
\def\eps{\varepsilon} \def\f{\varphi} \def\om{\omega} \def\Om{\Omega}
\def\gm{\gamma} \def\gep{\gm_\eps} \def\lm{\lambda} \def\lep{\lm_\eps}
\def\dl{\delta} \def\rb{\overline{r}} \def\rt{\tilde{r}} \def\al{\alpha}
\def\fh{\hat{f}} \def\yb{\overline{y}} \def\oet{\L(1-e^{-\mu\tau}\R)}
\def\mb{\overline{m}} \def\Mb{\overline{M}} \def\Mh{\hat{M}}
\def\bx{\mathbf{x}} \def\by{\mathbf{y}} \def\bS{\mathbf{S}}
\def\I{{\cal I}} \def\A{{\cal A}} \def\D{{\cal D}}\def\bc{\mathbf{c}}
\def\eq{equation} \def\de{differential \eq} \def\pde{partial \de}
\def\sol{solution} \def\pb{problem} \def\bdy{boundary} \def\fn{function}
\def\dde{delay \de} \def\ev{eigenvalue} \def\asy{asymptotic}
\def\bx{\mathbf{x}} \def\by{\mathbf{y}} \def\bS{\mathbf{S}}
\def\H{{\mathcal H}} \def\U{{\cal U}} \def\D{{\cal D}}\def\bc{\mathbf{c}}
\def\eq{equation} \def\de{differential \eq} \def\pde{partial \de}
\def\sol{solution} \def\pb{problem} \def\bdy{boundary} \def\fn{function}
\def\dde{delay \de} \def\ev{eigenvalue}
\def\R{\mathbf R}
\def\C{\mathbf C}

\author{
{\bf Dang Vu Giang}\\
Hanoi Institute of Mathematics\\
18 Hoang Quoc Viet, 10307 Hanoi, Vietnam\\
{\footnotesize          e-mail: $\<$dangvugiang@yahoo.com$\>$}\\
\1\\
}
\maketitle
\noindent {\bf Abstract.}   We are interested in computing  the spectral measure of Laplacean   operators  in Paley-Wiener space, the Hilbert space of all square integrable functions having Fourier transforms  supported in a compact set $K$, the closure of an open bounded set in $\R^N$. I is well-known that every differential operator is bounded in this space. Among others, we will prove that the spectrum of Laplace operator is the set
$$\{-|x|^2: x\in K\}.$$



\eject

\sect{Introduction }

\bigskip
\noindent
Let $K_0$ be a bounded open set in $\R^N$.
Let $K$ be the closure of $K_0$.
Then $K$ is a compact set.
Let $\mathcal H$ be the Hilbert space of all $f\in L^2(\R^N)$ such that the Fourier transform $\hat f$ of $f$ is supported in $K$. Let $\Delta$ denote the Laplace operator in $\R^N$. Then $\Delta$ generates a contraction $C_0$-semigroup in $L^2(\R^N)$ and
$$e^{t\Delta}f(x)=\frac1{(4\pi t)^{N/2}}\int_{\R^N}\exp\Bigl(-\frac{|x-y|^2}{4t}\Bigr)f(y)dy,$$
for every $f\in L^2(\R^N)$. 
(See [2,3,4].)
On the other hand, $\Delta$ is bounded and self-adjoint operator on $\mathcal H$, so
$$e^{t\Delta}f(x)=\sum_{k=0}^\infty\frac{t^k}{k!}\Delta^k f(x),$$
for every $f\in\mathcal H $. Let $\rho(\Delta)$ denote the resolvent set of $\Delta$ over $\mathcal H$ and let $R(\lambda,\Delta):=\bigl(\lambda-\Delta\bigr)^{-1}$ denote the resolvent operator.
By Hille-Yosida theorem, $(0,\infty)\subset\rho(\Delta)$ and
$$||R(\lambda,\Delta)||
\leqslant\frac1\lambda\qquad\hbox{ for all }\quad\lambda>0.$$
Moreover,
\EQ{resolvent}{\ARR{cll}{
R(\lambda,\Delta)
&\ds=\frac1\lambda\sum_{k=0}^\infty\frac{\Delta^k}{\lambda^k}
\\
&\ds=\sum_{k=0}^\infty(\mu-\lambda)^kR(\mu,\Delta)^{k+1}\quad\hbox{ (for }  0<\lambda<\mu)
\\   
&\ds=\int_0^\infty e^{-\lambda t} e^{t\Delta}dt.}}
Let
$$\sigma(\Delta):=\C\setminus\rho(\Delta)$$
denote the spectrum  of $\Delta$ over $\mathcal H$. Then, by a well-known theorem of Gelfand, $\sigma(\Delta)$ is non-empty compact set of the complex plane $\C$. Let $$r(\Delta):=\sup\bigl\{|z|:z\in\sigma(\Delta)\bigr\}$$
denote the spectral radius of $\Delta$ over $\mathcal H$. Then
$$r(\Delta)=\lim_{n\to\infty}||\Delta^n||^{1/n}.$$
Let
$$\hat f(x)=\frac1{(2\pi)^{N/2}}\int_{\R^N}e^{-ix\cdot y}f(y)dy$$
denote the Fourier transform of $f\in L^1(\R^N)$.
Classical Fourier Analysis will give that
$$\widehat{\Delta f}(x)=-|x|^2\hat f(x),$$
and consequently,
$$||\Delta||=\sup_{x\in K}|x|^2=r(\Delta)=:R.$$

\bigskip
\noindent
{\bf Theorem 1. }  {\it We have}
$$\sigma(\Delta)=\{-|x|^2: x\in K\}.$$

\bigskip
\noindent{\it Proof: } First, $\Delta$ is self-adjoint, so $\sigma(\Delta)\subset\R$. But the resolvent set of $\Delta$ is containing the positive semi-axis, so $\sigma(\Delta)\subset (-\infty,0]$.  On the other hand, the spectral radius of $\Delta$ is $R$, so $\sigma(\Delta)\subseteq[-R,0].$  But $\Delta$ has no eigenvector in $\H$,
so $\Delta$ has approximate point spectrum only. This means that $-\lambda\in\sigma(\Delta)$ if and only if there is a sequence $f_1,f_2,\cdots$ of unit vectors of $\H$ such that
$$\lim_{n\to\infty}||\Delta f_n+\lambda f_n||=0$$
(see \cite{arendt} page 64, for details).
We will prove that if $\lambda=|x|^2$ for some $x\in K$, there is a sequence $f_1,f_2,\cdots$ of unit vectors of $\H$ such that
$$\lim_{n\to\infty}||\Delta f_n+\lambda f_n||=0.$$
Indeed, let $x_0$ be an interior point of $K$ and $\lambda=|x_0|^2$. For $n$ large enough, the ball $B_n$ centered at $x_0$ with radius $1/n$ is contained in $K$.  Let $\varphi_n=:\hat{f_n}$ be a unit vector of  $L^2(\R^N)$ supported in  $B_n$ which is constant in this ball.  Then, $f_n$ is a unit vector of $\H$, and
\[\begin{aligned}
||\Delta f_n+\lambda f_n||^2=||\widehat{\Delta f_n}+\lambda\hat f_n||^2
&=\int_{B_n}\bigl(|x_0|^2-|x|^2\bigr)^2|\varphi_n(x)|^2 dx\\
&\leqslant\frac 1{n^2}\cdot\Bigl(2|x_0|+\frac1n\Bigr)^2,
\end{aligned}\]
which is certainly tends to 0 as $n\to\infty$. Similar proof will hold, if $x_0$ is a boundary point. Therefore,
$$\{-|x|^2:  x\in K\}\subseteq\sigma(\Delta).$$
Conversely, let $-\lambda\in\sigma(\Delta)$. Let
$f_0,f_1,\cdots$ of unit vectors of $\H$ such that
$$\lim_{n\to\infty}||\Delta f_n+\lambda f_n||=0.$$
Let $\varphi_n=\hat f_n$.  We have
$$||\Delta f_n+\lambda f_n||^2=||\widehat{\Delta f_n}+\lambda\hat f_n||^2
=\int_{K}\bigl(\lambda-|x|^2\bigr)^2|\varphi_n(x)|^2 dx
\ge\delta,$$
where $\delta=\min_{x\in K} \bigl(\lambda-|x|^2\bigr)^2= (\lambda-|x_0|^2)^2$  for some $x_0\in K$. But
$||\Delta f_n+\lambda f_n||\to 0$ as $n\to\infty$, so $\delta=0$, and consequently, $\lambda=|x_0|^2$ for some $x_0\in K$. The proof is now complete.

\bigskip
\noindent{\bf Remark. } The assumption that $K$ is the closure of an open set is very essential.
 If $K$ has an isolated point or a cluster point, the assertion of this Theorem does not hold!

\bigskip

\sect{Moments of Spectral Measures}

\bigskip
\noindent
Let $H_0=-\Delta$. Then $H_0$ is bounded, positive and self-adjoint operator of $\H$, because
$$\langle\Delta f, f\rangle=-\int_{\R^N}\bigl|\nabla f(x)\bigr|^2dx\leqslant 0\quad\hbox{ for all } f\in\H.
$$
Consider the spectral decomposition of $\Delta$:
$$\Delta=\int_{\sigma(\Delta)}\lambda dE_\lambda.$$
Let $\mu$ denote the spectral measure of $\Delta$ with respect with some function $v\in\mathcal H$ with $||v||=1$.
Then $\mu$ is probability measure supported in $[-R,0]$ and
$$\langle R(\lambda,\Delta)v,v\rangle=\int_{\sigma(\Delta)}\frac{d\mu(x)}{\lambda-x}$$ and
$$\langle e^{it\Delta}v,v\rangle=\int_{\sigma(\Delta)} e^{itx}d\mu(x).$$
Let $\{p_0,p_1,p_2,\cdots\}$ be the system of orthonormal polynomials with respect to $\mu$. Then there are two sequences $\{a_0,a_1,a_2\cdots\}$ and $\{b_0,b_1,b_2\cdots\}$ such that
$$xp_n(x)=a_{n-1}p_{n-1}(x)+b_np_n(x)+a_np_{n+1}(x)\qquad\hbox{ for all }
n=1,2,\cdots.$$
Let $A$ be a tridiagonal Jacobi matrix with the main diagonal $\{b_0,b_1,b_2\cdots\}$  and two other diagonals
$\{a_0,a_1,a_2\cdots\}$:

\begin{equation*}
A=\left(
\begin{matrix}
b_0&a_0&0&0&\cdots\\
  a_0&b_1&a_1&0&\cdots\\   
   0&a_1&b_2&a_2&\cdots \\
0&0&a_2&b_3& \cdots\\
\vdots& \vdots& \vdots& \vdots& \ddots
  \end{matrix}
\right)
\end{equation*}   
Let $\underline e_0=(1,0,0,\cdots)$ be the first vector in $\ell^2(\mathbf N_0)$. Then $A$ as a bounded linear operator in $\ell^2(\mathbf N_0)$ will have $\mu$ as spectral measure with respect to $\underline e_0$, i.e.,
$$\langle e^{itA}\underline e_0,\underline e_0\rangle=\int_{\sigma(\Delta)} e^{itx}d\mu(x).$$
Hence, $\Delta$ and $A$ are spectrally equivalent. Since $\Delta$ has no eigenvalue, so $$\sum_{n = 0}^{\infty}\mid p_n(x)\mid^2 = \infty \qquad \mbox{for all } x \in \C.$$ It is well known that the measure $\mu$ is absolutely continuous, i.e.  $d\mu(\xi) = \varphi(\xi)d\xi$, where $\varphi(\xi)> 0$ on $\sigma(\Delta)$.
If $K$ is the closed unit ball of $\R^N$ and
$$v(x)=\frac1{(2\pi)^{N/2}\sqrt{|K|}}\int_{K}e^{-ix\cdot y}dy,$$
then $\sigma(\Delta)=[-1,0]$ and
$$\int_{-1}^0\xi^k\varphi(\xi)d\xi = (-1)^k\frac{N}{N + 2k}, \quad\mbox{for } k = 0, 1, 2, \cdots.$$
Hence, $\varphi (\xi )=\frac{N}{2}{{\left| \xi  \right|}^{\frac{N}{2}-1}}$ is the density of  the spectral measure.

\bigskip
\par\noindent{\bf Acknowledgement.} 
Deepest appreciation is extended towards the NAFOSTED  (the National Foundation for Science and Techology Development in Vietnam) for the financial support.


\bigskip

\bigskip

\begin{thebibliography}{99}
\bibitem{arendt} W. Arendt et al, "One-parameter semi-groups of operators",  Lecture
Notes in Mathematics 1184,  Springer-Verlag (1986).
\bibitem{Aupetit} B. Aupetit, "A primer on spectral theory", Springer-Verlag (1991).
\bibitem{berthier} A.M. Berthier, "Spectral theory and wave operators for the Schr\"odinger equation", Pitman Advanced Publishing Program (1982).
\bibitem{goldstein} J.A. Goldstein, "Semigroups of Operators and Applications", Oxford University Press (1985).


\end{thebibliography}
\end{document}